\documentclass[10pt]{amsart}
\usepackage[all]{xy}
\usepackage{mathrsfs}
\usepackage{amssymb}
\usepackage{array}

\usepackage{color, graphicx}
\usepackage[all]{xy}

\usepackage{calrsfs}
\usepackage{times}
\usepackage[scaled=0.92]{helvet}

\newcommand{\Spe}{{\rm Spec\,}}

\newcommand{\A}{{ \rm Aut }}

\newcommand{\F}{{\mathbb{F}}}

\newcommand{\Z}{{\mathbb{Z}}}

\renewcommand{\mod}{{\;\rm mod\;}}

\def\Aut{\operatorname{Aut}}
\def\GL{\operatorname{GL}}
\def\gen#1{\langle #1\rangle}

\newtheorem{pro}{Proposition}[section]
\newtheorem{lemma}[pro]{Lemma}
\newtheorem{cor}[pro]{Corollary}
\newtheorem{theorem}[pro]{Theorem}

\theoremstyle{definition}
\newtheorem{rem}[pro]{Remark}
\newtheorem{orism}[pro]{Definition}

\begin{document}
\title{{Automorphisms of hyperelliptic modular curves $X_0(N)$ in positive 
characteristic}}

\author{
Aristides Kontogeorgis} 
\address{Max-Planck-Institut f\"ur Mathematik Vivatsgasse 7
 53111 Bonn,  and  \\Dept of Mathematics, Univ. of the Aegean, 83200, Samos, Greece.
}
\email{kontogar@aegean.gr}
\author{Yifan Yang}
\address{Max-Planck-Institut f\"ur Mathematik Vivatsgasse 7
 53111 Bonn,  and Department of Applied Mathematics \\
  National Chiao Tung University \\
  Hsinchu 300 \\
  TAIWAN}
\email{yfyang@math.nctu.edu.tw}
\date\today
\begin{abstract} 
  We study the automorphism groups of the reduction $X_0(N) \times
  \bar{\mathbb{F}}_p$ of a modular curve $X_0(N)$ over primes 
 $ p\nmid N$.
\end{abstract}

\subjclass[2000]{Primary 14H37; secondary 11F11, 11G20}
\maketitle
\section{Introduction}
Let $\mathcal{X} \rightarrow S$ be a family of curves of genus $g\ge
2$ over a base scheme $S$. For every point $P:\Spe k \rightarrow S$,
we will consider the {\em absolute} automorphism group of the fibre
$P$ to be the automorphism group
$
 \A_{\bar{k}}( \mathcal{X} \times_S \Spe \bar{k})
$
where $\bar{k}$ is the algebraic closure of $k$. Any automorphism
$\sigma$ acts like the identity on $\bar{k}$, so in our setting there
is no $\mathrm{Gal}(\bar{k}/k)$ contribution to the automorphism group
of any special fibre. The following theorem due to P. Deligne,
D. Mumford \cite[lemma I.12]{DelMum} compares the automorphism groups
of the generic and special fibres:
\begin{pro}[\cite{DelMum}]
  Consider a stable curve $\mathcal{X}  \rightarrow S$ over a scheme
  $S$ and let $\mathcal{X}_\eta$ denote its generic fibre. Every
  automorphism $\phi:\mathcal{X}_\eta \rightarrow \mathcal{X}_\eta$ can
  be extended to an automorphism $\phi:\mathcal{X} \rightarrow
\mathcal{X}$.
\end{pro}
Of course the special fibre of $\mathcal{X}$ might possess
automorphisms that can not be lifted to the generic fibre. For example
the Fermat curve
\[
x^{p^s+1}+y^{p^s+1}+z^{p^s+1}=0,
\] 
can be considered as a stable curve over $\Spe \Z[\frac{1}{p^s+1}]$
and has automorphism group $(\mu_n \times \mu_n) \rtimes S_3$ for all
geometric fibres above the primes $q\neq p$ but
$\mathrm{PGU}(3,p^{2s})$ for the prime $p$ \cite{Tze:95},
\cite{Leopoldt:96}. A special fibre $\mathcal{X}_p:=\mathcal{X}
\times_S S/p$ with $\A(\mathcal{X}_p)> \A(\mathcal{X}_\eta)$ will be
called \emph{exceptional}. In general we know that there are finite
many exceptional fibres and it is an interesting problem to determine
exactly the exceptional fibres.

There are some results towards this problem for some curves of
arithmetic interest. A. Adler \cite{Adler97} and C.S. Rajan
\cite{Rajan98} proved for the modular curves $X(N)$,  that
$X(11)_3:=X(11)\times_{ \Spe \Z} \Spe \mathbb{F}_3$ has the Mathieu
group $\mathrm{M}_{11}$ as the full automorphism group.
C. Ritzenthaler in \cite{Ritzenthaler03} and P. Bending, A. Carmina,
R. Guralnick \cite{BeCaGu} studied the automorphism groups of the
reductions $X(q)_p$ of modular curves $X(q)$ for various primes $p$.
It turns out that the reduction $X(7)_3$ of $X(7)$ at the prime $3$
has an automorphism group $\mathrm{PGU}(3,3)$, and $X(7)_3$ and $X(11)_3$
are the only cases where $\A X(q)_p > \A X(q)\cong \mathrm{PSL}(2,p)$.

In this paper we will investigate some modular curves of the 
form $X_0(N)$. Igusa \cite{Igusa59} proved that $X_0(N)$ has a
non-singular projective model which is defined by equations over
$\mathbb{Q}$ whose reduction modulo primes $p, p\nmid N$ are also
non-singular, or in a more abstract language that there is a proper
smooth curve $\mathcal{X}_0(N) \rightarrow \Z[1/N]$ so that for $p \in
\Spe \Z[1/N]$ the reduction
$\mathcal{X}_0(N) \times_{\Spe \Z} \mathbb{{F}}_p$ is the moduli space
of elliptic curves with a fixed cyclic subgroup of order $N$.

The automorphism group of the curve $X_0(N)$ at the generic fibre is
now well understood. Let us call an automorphism $\sigma$ of $X_0(N)$
\emph{modular} if $\sigma$ arises from the normalizer of $\Gamma_0(N)$
in ${\rm PGL}^+_2(\mathbb{Q})$. Then the results of
\cite{Elkies90,Kenku-Momose,Ogg:77} can be summarized as follows.

\begin{theorem}[\cite{Elkies90,Kenku-Momose,Ogg:77}] Assume that
  $X_0(N)$ has genus $g\ge 2$. Then all the automorphisms of $X_0(N)$
  are modular, with only two exceptions, namely the
  cases $N=37,63$, which have extra non-modular involutions.
\end{theorem}

\begin{rem} Note that if $4\nmid N$ and $9\nmid N$, then all modular
automorphisms are of the Atkin-Lehner type
\cite[Theorem 8]{Atkin-Lehner}. In particular, if $4\nmid N$ and
$9\nmid N$, then the automorphism group is an elementary abelian
$2$-group.
\end{rem}

Denote by $A(N,0)$ the absolute automorphism group of
$\mathcal{X}_0(N)$ at the generic fibre and by $A(N,p)$ the absolute
automorphism group at the reduction at the prime $p$. The problem of
determining the exact primes for which $A(N,p)> A(N,0)$ seems
difficult in general. However, when $X_0(N)$ is hyperelliptic, the
situation is relatively simple because the function field of a
hyperelliptic curve has a unique genus zero subfield of degree $2$,
and the problem of determining the automorphism group of a
hyperelliptic curve is essentially the same as that of determining the
automorphisms of a projective line that permute a set of marked
points. Therefore, as a starting point for our general problem of
studying automorphisms of $X_0(N)$ in positive characteristics, in this
note, we will investigate the automorphism groups of hyperelliptic
modular curves $X_0(N)$.

\begin{theorem} Table \ref{Table: Results} is the complete list of
  integers $N$ and primes $p$ such that the reduction of the
  hyperelliptic modular curve $X_0(N)$ modulo $p$ has exceptional
  automorphisms.
\end{theorem}

Here in the table, the notation $D_n$ denotes the dihedral group of
order $2n$,
\[
A:=\gen{ a,b,c \mid
 c^2,
    b a^{-2}  b^{-1} a^{-1},
    b^{-1}  a^3  b  a^{-1},    
b  a^{-1} c  b^{-1}  a^{-1}  c  a^{-1} c,
    \left(a^{-1}  b^{-1}  c  b^{-1}\right)^2
},
\] 
is a group of order $672$,
$$
 B:=\gen{ a,b,c\mid  c^2,
    a^{-5},
    b^{-1}  a^{-2}  b  a, 
    (c  b^{-1})^3, 
    a^{-1}  b  c  a^2  c  a  c}.
$$
is a group of order $240$, and
$$
  V_n=\gen{x,y \mid x^4,y^n,(xy)^2,(x^{-1}y)^2}
$$
is a group of order $4n$. Moreover, in the case of $X_0(37)$ in
characteristic $2$, the notation $E_{32-}$ represents the extraspecial
group $E_{32-}=(D_4\times Q_8)/\gen{(a,b)}$, where $a$ and $b$ denote
the nontrivial elements of the centers of $D_4$ and $Q_8$,
respectively.

\begin{table} \caption{Automorphism groups of $X_0(N)$ in positive
    characteristics}
\label{Table: Results}
\extrarowheight3pt
\begin{tabular}{cccp{3.3cm}p{2.8cm}} \hline\hline
  $N$ & Genus & Generic $\Aut$ & Exceptional primes & Exceptional $\Aut$
    \\ \hline\hline
  $22$ & $2$ & $(\Z/2\Z)^2$ & $3,29$ \newline $101$ &
  $D_6$ \newline $D_4$ \\ \hline
  $23$ & $2$ & $\Z/2\Z$ & $3,13,29,43,101,5623$ & $(\Z/2\Z)^2$ \\ \hline
  $26$ & $2$ & $(\Z/2\Z)^2$ & $7,31$ \newline $41,89$ 
    & $D_6$ \newline $D_4$ \\ \hline
  $28$ & $2$ & $D_6$ & $3$ \newline $5$ \newline $11$
    & $\GL(2,3)$
       \newline $B$ \newline $V_6$ \\ \hline
  $29$ & $2$ & $\Z/2\Z$ & $19$ \newline $5,
       67,137,51241$ & $D_4$ \newline $(\Z/2\Z)^2$ \\ \hline
  $30$ & $3$ & $(\Z/2\Z)^3$ & $23$ & $V_8$ \\ \hline
  $31$ & $2$ & $\Z/2\Z$
    & $3$ \newline $5,11,37,67,131,149$
    & $D_4$ \newline $(\Z/2\Z)^2$ \\ \hline
  $33$ & $3$ & $(\Z/2\Z)^2$ & $2$ \newline $19$ \newline $47$
    & $\GL(2,2)\times(\Z/2\Z)$ \newline $(\Z/2\Z)\times(\Z/4\Z)$
      \newline $(\Z/2\Z)^3$ \\ \hline
  $35$ & $3$ & $(\Z/2\Z)^2$ & --- & --- \\ \hline
  $37$ & $2$ & $(\Z/2\Z)^2$
    & $2$ 
      \newline $7,31$ \newline $29,61$
    & $E_{32-}\rtimes(\Z/5\Z)$ 
      \newline $D_6$ \newline $D_4$ \\ \hline
  $39$ & $3$ & $(\Z/2\Z)^2$ & $5$ \newline $29$
    & $(\Z/2\Z)^3$ \newline $(\Z/2\Z)\times(\Z/4\Z)$ \\ \hline
  $40$ & $3$ & $(\Z/2\Z)\times D_4$ & $3$ & $V_8$ \\ \hline
  $41$ & $3$ & $\Z/2\Z$ & $17$ & $(\Z/2\Z)^2$ \\ \hline
  $46$ & $5$ & $(\Z/2\Z)^2$ & $3$ & $(\Z/2\Z)\times(\Z/4\Z)$ \\ \hline
  $47$ & $4$ & $\Z/2\Z$ & --- & --- \\ \hline
  $48$ & $3$ & $(\Z/2\Z)\times S_4$ & $7$ & $A$, $|A|=672$ \\ \hline
  $50$ & $2$ & $(\Z/2\Z)^2$ & $3$ \newline $37$
    & $D_6$ \newline $D_4$ \\ \hline
  $59$ & $5$ & $\Z/2\Z$ & --- & --- \\ \hline
  $71$ & $6$ & $\Z/2\Z$ & --- & --- \\ \hline\hline
\end{tabular}
\end{table}


{\bf Acknowledgments:} This paper was completed during the authors'
visit at Max-Planck Institut f\"ur Mathematik in Bonn. The authors
would like to thank the institute for its support and hospitality.

%
%
%
%
%
%

%
%

\section{Automorphisms in characteristic $p\neq 2$}

\label{sec2}
 According to \cite{Ogg:74}, there are exactly $19$ values
of $N$ such that $X_0(N)$ is hyperelliptic. The equations of the form
$y^2=f(x)$ for hyperelliptic $X_0(N)$ have been  computed by
several authors \cite{Galb:96,Gonzalez,ShimMa95}. They are
tabulated in Table \ref{hyplist}, along with their modular
automorphisms. Here, for a divisor $e$ of $N$ with $(e,N/e)=1$, we let
$w_d$ be the Atkin-Lehner involution corresponding to the normalizer
$$
  \begin{pmatrix}ae & b\\ cN & de\end{pmatrix}, \quad
  ade-bcN/e=1,
$$
of $\Gamma_0(N)$. For $X_0(28)$ and $X_0(40)$, the notation $w_{1/2}$
represents the automorphism coming from the normalizer
$$
  \begin{pmatrix}1&1/2\\0&1\end{pmatrix}.
$$
The notation $w_{1/4}$ in the case $X_0(48)$ carries a similar meaning.

\begin{table}
\caption{List of hyperelliptic curves $X_0(N)$ 
\label{hyplist}}
$$ \extrarowheight3pt
{\tiny
\begin{tabular}{|c||p{11cm}|} \hline\hline

$N$ & Equation/Automorphisms \\ \hline\hline

$22$ & $y^2=(x^3+4x^2+8x+4)(x^3+8x^2+16x+16)$ \newline
     $w_2:(x,y)\mapsto(4/x,8y/x^3), \quad w_{11}:(x,y)\mapsto(x,-y)$
     \\ \hline
$23$ & $y^2=(x^3-x+1)(x^3-8x^2+3x-7)$ \newline
     $w_{23}:(x,y)\mapsto(x,-y)$ \\ \hline
$26$ & $y^2=x^6-8x^5+8x^4-18x^3+8x^2-8x+1$ \newline
     $w_{13}:(x,y)\mapsto(1/x,y/x^3), \quad w_{26}:(x,y)\mapsto(x,-y)$
     \\ \hline
$28$ & $y^2=(x^2+7)(x^2+x+2)(x^2-x+2)$ \newline
     $w_4:(x,y)\mapsto((x+3)/(x-1),8y/(x-1)^3), \quad
      w_7:(x,y)\mapsto(x,-y), \quad w_{1/2}:(x,y)\mapsto(-x,-y)$
     \\ \hline
$29$ & $y^2=x^6-4x^5-12x^4+2x^3+8x^2+8x-7$ \newline
     $w_{29}:(x,y)\mapsto(x,-y)$
     \\ \hline
$30$ & $y^2=(x^2+4x-1)(x^2+x-1)(x^4+x^3+2x^2-x+1)$ \newline
     $w_2:(x,y)\mapsto((x+1)/(x-1),-4y/(x-1)^4), \quad
      w_3:(x,y)\mapsto(-1/x,-y/x^4)$, \newline
     $w_{15}:(x,y)\mapsto(x,-y)$
     \\ \hline
$31$ & $y^2=(x^3-6x^2-5x-1)(x^3-2x^2-x+3)$ \newline
     $w_{31}:(x,y)\mapsto(x,-y)$
     \\ \hline
$33$ & $y^2=(x^2+x+3)(x^6+7x^5+28x^4+59x^3+84x^2+63x+27)$ \newline
     $w_3:(x,y)\mapsto(3/x,-9y/x^4), \quad
      w_{11}:(x,y)\mapsto(x,-y)$
     \\ \hline
$35$ & $y^2=(x^2+x-1)(x^6-5x^5-9x^3-5x-1)$ \newline
     $w_{7}:(x,y)\mapsto(-1/x,y/x^4), \quad w_{35}:(x,y)\mapsto(x,-y)$
     \\ \hline
$37$ & $y^2=x^6+14x^5+35x^4+48x^3+35x^2+14x+1$ \newline
     $w_{37}:(x,y)\mapsto (1/x,y/x^3 )$
     \\ \hline
$39$ & $y^2=(x^4-7x^3+11x^2-7x+1)(x^4+x^3-x^2+x+1)$ \newline
     $w_3:(x,y)\mapsto(1/x,y/x^4), \quad
      w_{39}:(x,y)\mapsto(x,-y)$
     \\ \hline
$40$ & $y^2=x^8+8x^6-2x^4+8x^2+1$ \newline
     $w_5:(x,y)\mapsto(-1/x,-y/x^4), \quad
      w_8:(x,y)\mapsto((1-x)/(1+x),-4y/(x+1)^4)$ \newline
     $w_{1/2}:(x,y)\mapsto(-x,y)$
     \\ \hline
$41$ & $y^2=x^8-4x^7-8x^6+10x^5+20x^4+8x^3-15x^2-20x-8$ \newline
     $w_{41}:(x,y)\mapsto(x,-y)$
     \\ \hline
$46$ & $y^2=(x^3+x^2+2x+1)(x^3+4x^2+4x+8)(x^6+5x^5+14x^4+25x^3+28x^2+20x+8)$
     \newline
     $w_{23}:(x,y)\mapsto(x,-y), \quad w_{46}:(x,y)\mapsto(2/x,8y/x^6)$
     \\ \hline
$47$ & $y^2=(x^5+4x^4+7x^3+8x^2+4x+1)(x^5-5x^3-20x^2-24x-19)$ \newline
     $w_{47}:(x,y)\to(x,-y)$
     \\ \hline
$48$ & $y^2=(x^4-2x^3+2x^2+2x+1)(x^4+2x^3+2x^2-2x+1)=x^8+14x^4+1$ \newline
     $w_{1/4}:(x,y)\to(ix,y), \quad w_3:(x,y)\to(-1/x,-y/x^4),$ \newline
     $w_{16}:(x,y)\to((1-x)/(1+x),-4y/(1+x)^4)$
     \\ \hline
$50$ & $y^2=x^6-4x^5-10x^3-4x+1$ \newline
     $w_2:(x,y)\mapsto(1/x,y/x^3), \quad
      w_{50}:(x,y)\mapsto(x,-y)$
     \\ \hline
$59$ & $y^2=(x^3+2x^2+1)(x^9+2x^8-4x^7-21x^6-44x^5-60x^4-61x^3-46x^2-24x-11)$
     \newline
     $w_{59}:(x,y)\mapsto(x,-y)$
     \\ \hline
$71$ & $y^2= (x^7-3x^6+2x^5+x^4-2x^3+2x^2-x+1)
           (x^7-7x^6+14x^5-11x^4+14x^3-14x^2-x-7) $ \newline
     $w_{71}:(x,y)\mapsto(x,-y)$
     \\ \hline
\end{tabular}
}
$$
\end{table}

These models are not the Igusa curves. They have a singularity 
at infinity but are smooth otherwise. 
The automorphism group of a non-singular curve is a 
birational invariant and one can work with singular
models of the corresponding function fields.
But one has to be careful if he wants to 
compute reductions out of the singular model. 
\begin{pro}
Let $y^2=f(x)$ where $f(x)\in \Z[x]$ is a polynomial with integer 
coefficient such that $\deg(f) \equiv 0 \mod 2$.

Suppose that $f(x)$ has no multiple roots and 
denote $\rho_1,\ldots,\rho_s \in \bar{\mathbb{Q}}$ the set of roots of $f$.  
Let $\Delta \in \Z$ be the discriminant of the polynomial $f(x)$.
Then, the relative curve  
\begin{equation} \label{curvedef}
\mathcal{Y}:=\frac{\Z\left[\Delta^{-1} \right][x,y]}{\langle y^2=f(x) \rangle} \rightarrow 
\mathbb{Z}\left[ \Delta^{-1} \right] 
\end{equation}
is a  smooth family. The set of roots $\{\rho_i\}_{i=1,\ldots,s}$
gives rise to horizontal divisors $\bar{\rho}_i$ that intersect 
the generic fibre of $\mathcal{Y}$ at $\rho_i$. 
The intersections of $\bar{\rho}_i$ with the fibre $\mathcal{Y}_p$ 
at the prime $p$ are roots of the polynomial  $f \mod p$.

Let $\mathcal{X} \rightarrow \mathbb{Z}\left[ (2N)^{-1}\right]$ be the 
Igusa family (with the $2$-fiber removed if $2 \nmid N$).
For every $p$ such that $(p,\Delta)=1$ the function fields 
of the curves $\mathcal{X}_p$ and $\mathcal{Y}_p$ have the 
same automorphism group.
\end{pro}
\begin{proof}
To prove that the map given in eq. (\ref{curvedef}) 
is smooth we will show that $\mathcal{Y}_p$ is smooth for every $p$
so that $p\nmid \Delta$. 
Set $F(x)=y^2-f(x)$. We compute 
\[
\frac{\partial F}{\partial y}=2y,\;\;\; \frac{\partial F}{\partial x}=
\frac{\partial f}{\partial x}. 
\] 
Therefore, a nonsmooth point appears at $y=0$  (since $p\neq 2$) 
and at a double 
root of $f(x) \mod p$, i.e., only at primes dividing the discriminant.

The Igusa family has good reduction at the prime $2$ if $2\nmid N$. 
The hyperelliptic models in Table \ref{hyplist} have always bad
reduction at the prime $2$. Therefore the reduction of the curve
$\mathcal{Y}_2$ is not related to the Igusa curve. In the $2\nmid N$
case the special fibre of the Igusa model at the prime $2$,  is a
nonsingular curve and the cover  $\mathcal{X}_0(N)_2 \rightarrow
\mathbb{P}^1_{\bar{\mathbb{F}}_2}$  is given in terms of an
Artin-Schreier extension. The automorphisms in characteristic $2$ will
be treated separately in Section \ref{section: char 2}.


The Igusa curve is the  normalization of  $\mathcal{Y}$  at infinity at 
every prime  in $\Spe \Z[\Delta^{-1}]$. 
It is known  that every point $P$ 
in the generic fibre $\mathcal{Y}_\eta$ defines a unique thickening, i.e.,
 a horizontal branch divisor intersecting the generic fibre at $P$.
Since $\deg(f)\equiv 0 \mod 2$ the place at infinity is not ramified at  the 
cover $X_0(N) \rightarrow \mathbb{P}^1$. Therefore the set of branch points 
of the cover $X_0(N)  \rightarrow \mathbb{P}^1$ is contained in the 
generic fibre  
$\mathcal{Y}_\eta$ of $\mathcal{Y}$. Moreover the thickenings 
$\tilde{\rho}_i$ of each point $\rho_i$ in 
$\mathcal{Y}_\eta \otimes_{\mathbb{Q}} \bar{\mathbb{Q}} \subset 
 \mathcal{X}\otimes_{\mathbb{Q}} \bar{\mathbb{Q}}$
do not intersect the singular locus of the projective closure of $\mathcal{Y}$.

The automorphism group of $\mathcal{X}_p$ is equal to the 
automorphism group of the corresponding 
function field and $\mathcal{X}_p, \mathcal{Y}_p$ share the 
same function field.
\end{proof}

The automorphism group of hyperelliptic curves is a well studied 
object \cite{BrandtStichtenoth}, \cite{GutSha},\cite{shaska03}.
There is the hyperelliptic involution $j$ so that 
$\langle j \rangle$ is a normal subgroup of the 
whole automorphism group, and moreover $X/\langle j \rangle$ is 
isomorphic to the projective line. 
The quotient $\bar{G}=\A(X)/\langle j \rangle$ is called the 
reduced group. The reduced group is a finite subgroup 
of the group of automorphisms of the projective line and these groups are 
given by proposition \ref{valmad}.
We introduce the following notation: For every $w\in \A(X)$ we will denote 
by $\bar{w}$ the image of $w$ in the reduced group $\bar{G}:=\A(X)/\langle j \rangle$.

Consider the cover $X \rightarrow X/\langle j \rangle\cong \mathbb{P}^1_k$.
Let $P_1,\ldots,P_s \in \mathbb{P}^1_k$ be the branched points
of the cover $X \rightarrow \mathbb{P}^1_k$.  
We will call the set $P_1,\ldots,P_s$ the branched hyperelliptic locus.

We know that the reduced group $\bar{G}$ induces a permutation action on 
them. If the points $P_1,\ldots, P_s$ are in general position 
then there is no finite subgroup of $\mathrm{PGL}(2,k)$ permuting
them and $\A(X)=\langle j \rangle$. The existence 
of additional automorphisms is a matter of special configurations 
of points in the configuration space of $s$ points in the 
projective line. 
The problem of determining the primes $p$, so that 
the reduction of a hyperelliptic curve at that primes 
has more automorphisms, is then reduced to the 
problem of determining the primes at which the branch hyperelliptic 
locus becomes more symmetric.

We now return  to the theory of  hyperelliptic curves $X_0(N)$.
If the number $N$ is composite or $N=37$  then the group of  automorphisms  
at the generic fibre is bigger than $\Z/2\Z$. This means that we have non trivial reduced group 
at the generic fibre. 
This situation gives us a lot of information for the location of 
fixed points of any extra  automorphism in the reduction modulo $p$ and on 
the possible extra automorphism groups.  

For instance the  case of hyperelliptic curves that admit an 
extra involution is studied by  T. Shaska and 
J. Gutierrez \cite{shaska03},\cite{GutSha}.
They introduced the  theory of {\em dihedral invariants},
a theory that  allows us to compute 
every possible extra automorphism for hyperelliptic curves 
 of genera $2,3$, 
and also gives us a lot of information in the bigger genus cases.

If the reduced automorphism group at  the generic fibre is trivial 
then we use a brute force method in order to compute any extra 
automorphism group. This is a  demanding computational  problem 
that needs several days 
of processing time.

Once the reduced group $\bar{G}$ of a hyperelliptic curve is determined
the group $G$ is given in terms of an extension of groups 
\[
 1 \rightarrow \langle j \rangle \rightarrow  G \rightarrow \bar{G} \rightarrow 1.
\]
For a cohomological approach to the structure of $G$ we refer to 
\cite{Ko:99}. It is known the  group structure of $G$ depends on whether 
the fixed points of $\bar{G}$ are in the branch locus of the cover
$X_0(N) \rightarrow X_0(N)^{\langle j \rangle}=\mathbb{P}^1_k$ 
\cite{BrandtStichtenoth}.

\begin{pro} \label{valmad}
Let $k$ be an algebraically closed field of characteristic $p \geq 0$.
Let $G$ be a  finite subgroup of $\A(\mathbb{P}^1_k)=\mathrm{PGL}(2,k)$, 
let $P_1,\ldots,P_r$ denote the number of branch points in the cover $\mathbb{P}^1_k
\rightarrow \mathbb{P}^1_k$. 
Denote the ramification degree of $Q \mapsto P_i$ by $e_i$. Then $G$ is 
one of the following groups:
\begin{enumerate}
\item Cyclic group $\Z/n\Z$ of order $n$ relatively prime to $p$ with $r=2$, 
and $e_1=e_2=n$.
\item Elementary abelian $p$-group with $r=1$, $e_1=|G|$.
\item Dihedral group $D_n$ of order $2n$, with $p=2$, $(p,n)=1$, $r=2$, 
      $e_1=2$, $e_2=n$, or $p\neq 2$, $(p,n)=1$, $r=3$, $e_1=e_2=2$, 
$e_3=n$.
\item Alternating group $A_4$ with $p\neq 2,3$, $r=3$, $e_1=2$, $e_2=e_3=3$
\item Symmetric group $S_4$ with $p\neq 2,3$, $r=3$, $e_1=2$, $e_2=3$, $e_3=4$.
\item Alternating group $A_5$ with $p=3$, $r=2$, $e_1=6$, $e_2=5$, or $p\neq 2,3,5$ $r=3$, $e_1=2$, $e_2=3$, $e_3=5$.
\item Semidirect product of an elementary abelian $p$-group of order $p^t$
with a cyclic group $\Z/n\Z$ of order $n$ with $n|(p^t-1)$, $r=2$, $e_1=|G|$, 
$e_2=n$.
\item $\mathrm{PSL}(2,p^t)$ with $p\neq 2$, $r=2$, $e_1=p^t(p^t-1)/2$, 
$e_2=(p^t+1)/2$.
\item $\mathrm{PGL}(2,p^t)$ with $r=2$ $e_1=p^t(p^t-1)$, $e_2=p^t+1$.
\end{enumerate}
\end{pro}
\begin{proof}
\cite[Th. 1]{Val-Mad:80}
\end{proof}

\begin{rem} Observe that the groups $A_4,S_4,A_5,\mathrm{PSL}(2,p^t),
\mathrm{PGL}(2,p^t)$ contain a dihedral group $D_n$, for $(n,p)=1$.
\end{rem}

Let $f(x)$ be a polynomial of degree $s$ with roots $\rho_1,\ldots,\rho_s$.
A simple computation shows that 
\begin{lemma}
Let $A=\begin{pmatrix} a & b \\ c & d \end{pmatrix}$ be 
an invertible matrix. 
If $c=0$ and  $\rho_i c \neq a$ for all roots of $f$ then:
\begin{equation} \label{11}
f\left(\frac{ax+b}{cx+d} \right)=\frac{c}{(cx+d)^s}f(a/c) 
\prod_{i=1}^s \left(x- \frac{\rho_id-b}{-\rho_ic+a}\right).
\end{equation}
If $c=0$ then $\rho_i c \neq a$ since $A$ is invertible and  
\begin{equation} \label{11a}
f\left(\frac{ax+b}{d} \right)=
\left(\frac{a}{d}\right)^s \prod_{i=1}^s \left(x- \frac{\rho_i d-b}{a} \right),
\end{equation}
If there is some $\rho_{i_0}$ such that $c \rho_i=a$ then $c\neq 0$ 
and 
\begin{equation} \label{11b}
f\left(\frac{ax+b}{cx+d} \right)=
\frac{(d \rho_{i_0} -b)}{(cx+d)^s}
 \prod_{i=1,i\neq i_0}^s  (a -\rho_i c)
\prod_{i=1,i\neq i_0}^s 
 \left(x- \frac{\rho_id-b}{-\rho_ic+a}\right).
\end{equation}
\end{lemma}

\begin{orism}
Let $A$ denote the M\"obius transformation 
$x \mapsto (ax+b)/(cx+d)$. We will denote by $f^*_A$ the 
polynomial 
\[
f^*_A(x):=\prod_{i=1, \rho_i c \neq a}^s \left(x- \frac{\rho_id-b}{-\rho_ic+a}\right).
\]
Notice that $\deg(f)=\deg(f^*_A)$ if $c \rho_i \neq a$ for all roots $\rho_i$ 
of $f$,  or $\deg(f)=\deg(f^*_A)+1$ if there is a root $\rho_i$ such that 
$c \rho_i=a$.
\end{orism}

\begin{lemma}
Let $A$ be a M\"obius transformation and set 
 $s=\deg f$. If $s\equiv 0 \mod 2$ then 
the curves $y^2=f(x)$ and $y^2= f^*_A(x)$ are isomorphic, over a quadratic extension 
of $\mathbb{Q}$.
\end{lemma}
\begin{proof}
The two curves become isomorphic under the change of variables 
\[(x,y) \mapsto 
\left(\frac{ax+b}{cx+d},\frac{y C^{1/2}}{(cx+d)^{s/2}  } \right),
\]
where $C$ is a constant depending on which of the cases in (\ref{11}),(\ref{11a}),(\ref{11b}) we are.
\end{proof}

If $f(x)=\sum_{\nu=0}^s a_\nu x^nu$ is a polynomial with $a_0\neq 0$ then 
we will denote by $f^*(x)$ the {\em reciprocal} polynomial given 
by $f^*(x)=a_0^{-1}x^s f(1/x)=a_0^{-1}\sum_{\nu=0}^s a_\nu x^{s-\nu}$.

If $2|s$ and $f(0)\neq 0$ then 
the hyperelliptic curves 
\[
y^2=f(x) \mbox{ and } y^2=f^*(x)
\]
are isomorphic.

\begin{lemma}
We will denote by $\sigma$ the automorphism 
of $\mathbb{P}^1_k$ sending  $x$ to $-x$. 
Consider the hyperelliptic curve 
\[
X: y^2=f(x)=\sum_{\nu=0}^s a_{2\nu} x^{2\nu},
\]
that has $\sigma$ in the reduced automorphism group.
There  is a cyclic group $\langle \sigma \rangle < C_d<\mathrm{PGL}(2,k)$
 that is a subgroup of the reduced automorphism group of $X$ 
 if and only if  $d|s$ and $a_\delta=0$ for all $d|\delta$. 
\end{lemma}

%

%
%
%
\section{Curves with non trivial reduced group}
%
%
\label{sec3}

\subsection{Curves of Genus $2$ with an extra involution.}
In this case we consider the set of curves $X_0(N)$ that are of genus $2$ and 
have reduced group isomorphic to $\Z/2\Z$.
This is the case for $N=22,26,28,37,50$. 

For genus two curves with reduced group $\Z/2\Z$ there is a well
developed theory due to Guttierez, Shaska \cite{GutSha},
\cite{shaska03}, namely the 
theory of {\em dihedral invariants} that reduce the computation 
of the tame part of the automorphism group to the computation of
several invariants of the curve.

Since $X_0(N)$ has reduced group $\Z/2\Z$, we can find a model 
of the curve so that the generator $\sigma$ of the 
reduced group acts like $x \mapsto -x$. Thus, the 
model of our curve is of the form:
\begin{equation} \label{eqdiinv}
 y^2=x^{2g+2} + a_1 x^{2g} +\cdots + a_g x^2 + 1.
\end{equation}
 The dihedral invariants are then given by $u_i:=a_1^{g-i+1}a_i+ a_g ^{g-i+1} a_{g-i+1}$ 
for $i=1,\ldots,g$. In particular for a genus two curve the dihedral invariants  are given by 
$u_1=a_1^3+a_2^3$, $u_2=2a_1a_3$.
Let $V_n$ denote the group 
\[
V_n:=\langle x,y \mid x^4,y^n,(xy)^2, (x^{-1} y)^2 \rangle.
\]
 In \cite[exam. 5.2]{ShaskaACM03}, \cite{ShaskaVolklein} T. Shaska and 
H. V\"olklein proved that 
the automorphism group is isomorphic to 
\begin{enumerate}
 \item \label{aa}
$V_6$ if and only if $(u_1,u_2)=(0,0)$ or $(u_1,u_2)=(6750,450)$
\item \label{bb} 
\begin{enumerate}
\item $\mathrm{GL(2,3)}$ if and only if $(u_1,u_2)=(-250,50)$ and $p\neq 5$
\item $B$ if and only if $(u_1,u_2)=(-250,50)$ and $p=5$
\end{enumerate}
\item $D_6$ if and only if $u_2^2-220 u_2-16 u_1+4500=0$,
\item $D_4$ if and only if $2u_1^2-u_2^3=0$ for $u_2\neq 2,18,0,50,450$. 
\end{enumerate}
(Cases $0,450,50$ are reduced to Cases (\ref{aa}) and (\ref{bb}) ).
The group $B$ mentioned above is given by:
\[
B:=\langle a,b,c\mid  c^2,
    a^{-5},
    b^{-1}  a^{-2}  b  a, 
    (c  b^{-1})^3, 
    a^{-1}  b  c  a^2  c  a  c \rangle.
\]

Using the exact action of the generator of the reduced group given in
Table  \ref{hyplist} 
we can find a model  of the form (\ref{eqdiinv}) by  diagonalizing the $2 \times 2$ matrix
representing the M\"obius transformation.
The dihedral invariants are then computed (\ref{tabdihinv2})
\begin{table} 
 \caption{Dihedral invariants of curves of genus 2 with extra involutions \label{tabdihinv2}}
$$ \extrarowheight2pt
\begin{array}{c|cc} \hline\hline
 N & u_1 & u_2 \\ \hline \hline
22 &  -17322/14641  & 130/121 \\
26 & -4351/2704 & 15/13 \\
28 & 43625/784 & 125/7 \\
37 &  -25642/1369 &-198/37  \\
50 & -135/16 & -5 \\ \hline\hline
\end{array}
$$
\end{table}
and this information allows us to compute the automorphism groups given in table \ref{Table: Results}.

%
%
%
\subsection{Curves of genus $3$ with an extra involution}
A similar approach to the genus $2$ in terms of the dihedral invariants can be given 
for the hyperelliptic curves with an involution in the reduced group that have genus $3$.
A complete list of all possible automorphism groups  that can appear in the genus $3$ case
together with necessary conditions on the dihedral invariants is given in 
the following theorem due to J.Gutierrez, D. Sevilla and T. Shaska
\cite{GutSevSha05}.
\begin{theorem}
Consider a hyperelliptic curve $X$ of genus $3$ with an 
extra involution. Let $G$ denote the full 
automorphism group $\A (X)$ of $X$  and $\bar{G}$ the reduced automorphism 
group. If the curve has $G,\bar{G}$ as in the first two columns of
table \ref{tableAut3} then the conditions given in  the 3rd column
of table \ref{tableAut3} are 
satisfied, 
where 
\begin{eqnarray*}
E_1(u_1,u_2) & :=&588 u_2-5(u_3-8)(9u_3-1024),\\
E_2(u_1,u_2) & := &7^3u_1-\frac{9}{8}u_3^3-\frac{873}{2}u_3^2+\frac{149504}{9}u_3-\frac{1048576}{9}
\end{eqnarray*}
and $V_8,U_6$ are the groups with presentations:
\begin{eqnarray*}
V_8 &:= &\langle x,y \mid x^4,y^8,(xy)^2,(x^{-1}y)^2 \rangle \\
U_6 &:= & \langle x,y \mid x^2,y^{12},xyxy^7 \rangle.
\end{eqnarray*}
\end{theorem}
\begin{table}
\caption{Possible Automorphism groups in the $g=3$ case. \label{tableAut3}}
{\tiny \extrarowheight3pt
\begin{tabular}{c|c|p{7.7cm}} \hline\hline
$G$ & $\bar G$ & Conditions \\ \hline\hline
$(\Z/2\Z)^3$ & $(\Z/2\Z)^2$ & $2u_1-u_3=0$\\ \hline
$(\Z/2\Z) \times D_8$ & $D_8$ & $2u_1-u_3=0,a_1=a_3$ \\ \hline
$V_8$ & $D_{16}$ & $2u_1-u_3=0,a_1=a_2=a_3=0$\\ \hline
$(\Z/2\Z) \times S_4$ & $ S_4$ & 
($2u_1-u_3=0,a_1=a_3$ or $ E_1(u1,u_2)=E_2(u_1,u_2)=0$) and
  \newline
$((u_1,u_2,u_3)=(0,196,0)$ or $(81 u_1,27u_2,8u_3)=(8192,-1280,128)$)
\\ \hline
$D_{12}$ & $D_6$ & $E_1(u_1,u_2)=E_2(u_1,u_2)=0$\\
\hline
$\Z/2\Z \times \Z/4\Z$ & $D_4$ & $ 2u_1+u_3^2=0$ \\
\hline
$U_{6}$ & $D_{12}$ & 
($2u_1+u_3^2=0$ or $E_1(u1,u_2)=E_2(u_1,u_2)=0$) and\newline
$u_2=0$. \\
\hline\hline
\end{tabular}
}
\end{table}

The curves $X_0(N)$ that are hyperelliptic of genus $3$ and have an involution in the 
reduced group correspond to $N\in \{39,40,48, 33,35,30\}$.
We compute first a hyperelliptic model of our curves so that the 
generator $\sigma$ of an extra involution is given by 
$\sigma:x \mapsto -x$. These models are given in table \ref{tabx2g3}.

\begin{table} \extrarowheight3pt
\caption{Models $f(x)$ for $N=33,35,39,40,48$. \label{tabx2g3} }
\begin{tabular}{l|l} \hline\hline
$N$ & $f(x)$ \\
\hline
$30$ & $  {x}^{8}+ \frac{\left( 276+184\,\sqrt {2} \right)}{\left( -540\,\sqrt {2}-765 \right)}
 {x}^{6}-46\,{x}^{4}+ 
\frac{\left( -184\,\sqrt {2}+276 \right)}{\left( -540\,\sqrt {2}-765 \right)}{x}^{
2}-\frac{765+540\,\sqrt {2}}{\left( -540\,\sqrt {2}-765 \right)}$ \\
$33$ & ${x}^{8}+{\frac { \left( -240\,\sqrt {3}+508 \right) {x}^{6}}{-264\,
\sqrt {3}+473}}+342\,{x}^{4}+{\frac { \left( 508+240\,\sqrt {3}
 \right) {x}^{2}}{-264\,\sqrt {3}+473}}+{\frac {473+264\,\sqrt {3}}{-
264\,\sqrt {3}+473}}$
\\
$35$ & $5\,{x}^{8}+ \left( 140+128\,i \right) {x}^{6}-34\,{x}^{4}+ \left( 140-
128\,i \right) {x}^{2}+5$
 \\
$39$ & $27\,{x}^{8}-2^2\cdot 97\,{x}^{6}+2\cdot29\,{x}^{4}+2^2\cdot 11\,{x}^{2}+3$\\
$40$ & $x^8-18x^4+1$\\
$48$ & $x^8+14x^4+1$\\
\hline\hline
\end{tabular}
\end{table} 
The dihedral invariants
in the case of genus $3$ curves  with an extra involution 
are given by $u_1=a_1^4+a_3^4$, $u_2=(a_1^2+a_3^2)a_2$, $u_3=2a_1a_3$, 
where 
\[y^2=x^8+a_1x^6+a_2 x^4+a_3x^2+1\]
is a normalized model of the hyperelliptic curve.
The dihedral invariants for the hyperelliptic curves $X_0(N)$ of genus $3$ with 
extra involution are given in table \ref{dig33}.
 \begin{table}
\caption{ Exceptional primes  for $N=30,33,35,39,40,48$. \label{dig33} }
{ \extrarowheight2pt
\begin{tabular}{c|p{1.8cm}p{1.8cm}p{2.5cm}p{1.8cm}}
\hline\hline
$N$ &  Factors of\newline $2u_1 +u_g^{2}$ & 
Factors of \newline $2u_1 -u_g^{2}$
 &  Possible excep-\newline tional primes  &
Exceptional primes  \\
\hline
$30$ & $2,23,17$ & $2,17$ & $23,17$ & $17$\\
$33$ & $2,19,31,103$ & $2,3,47$ & $19,31,47,103$  & $19,47$\\
$35$ & $3,67$ & $2,7$           & $3,67$ &  \\
$39$ & $2,29,181 $ & $2,5,13 $  & $5,29,181$ & $5,29$ \\
$40$ & $0$             & $0$          & $3$ & $3$ \\
$48$ & $0$              & $0$      & $7$  & $7$  \\
\hline\hline
\end{tabular}
}
\end{table}

For $N=39$ we see that possible exceptional primes are $5,29,181$ 
and then by reducing the coefficients modulo each of these primes 
we see that $A(39,5)\cong(\Z/2\Z)^3$, $A(39,29)\cong(\Z/2\Z) \times(\Z/4\Z)$.
The prime $181$ is not exceptional. 
For  $N=40$ the possible exceptional primes are $3,7$ and 
 $A(40,0)=(\Z/2\Z) \times D_8$ and   
 $A(40,3)\cong V_8$ while for $A(40,7)\cong (\Z/2\Z) \times S_4$.
For  $N=48$ the possible exceptional primes is $7$ and 
we have $A(48,0)=(\Z/2\Z) \times S_4$, $A(48,7)=V_8$.
The full group of  $\mathcal{X}_0(48)_7$
was studied by the first author in \cite{Ko:98} and is isomorphic 
to an extension of $\mathrm{PGL}(2,7)$ by $\Z/2\Z$.
Using magma we compute that this group admits the 
following presentation:
\[
A:=\left\langle a,b,c \mid
 c^2,
    b a^{-2}  b^{-1} a^{-1},
    b^{-1}  a^3  b  a^{-1},    
b  a^{-1} c  b^{-1}  a^{-1}  c  a^{-1} c,
    \left(a^{-1}  b^{-1}  c  b^{-1}\right)^2
\right\rangle.
\]

For $N=30,33,35$ the situation is a little more difficult, since 
the normal model is not defined over $\Z$ but in the 
principal ideal domains $\Z[\sqrt{2}],\Z[\sqrt{3}],\Z[i]$
respectively.
We compute that 
the dihedral invariants in the cases $N=33,35$ are in $\mathbb{Q}$ and in  
table \ref{dig33} we present the prime factors of the numerator.
For the $N=30$ case the dihedral invariants are in $\mathbb{Q}[\sqrt{2}]$,
and we compute that the principal ideal generated by the 
numerators of $2u_1^2-u_g^2$ and $2u_1^2+u_g^2$ have  the 
prime ideals  $I_2,I_{23,1},I_{23,2}$ and $I_2, I_{23,1},I_{23,2},
I_{17,1},I_{17,2}$, 
where $I_2=\langle 2,\sqrt{2}\rangle_{\Z[\sqrt{2}]}$ and 
$I_{23,i},I_{17,i}$ are the prime ideals that extend the prime 
ideals  $23\Z,17\Z$ of the ring of rational integers.
Since both conjugate prime ideals $I_{23,i}$, $i=1,2$ (resp. 
$I_{17,i}$) are divisors of $2u_1^2-u_g^2$ (resp. $2u_1^2+u_g^2$)
we see that $17,23|(2u_1+u_g^2)$ and $17|(2u_1+u_g)$.
The possible exceptional primes are $p=17,33$. 
We reduce the 
coefficients modulo  $p=17,33$ and we compute the 
automorphism group. 
It turns out that $A(30,17)\cong V_8$, $A(33,19)\cong(\Z/2\Z) \times (\Z/4\Z)$,
$A(33,47) \cong (\Z/2\Z)^3$.

\subsection{Elementary abelian groups}
%
%
We will study now whether there might be an elementary 
abelian subgroup $E\cong(\Z/p\Z)^s$ as a subgroup of the reduced group
of a hyperelliptic curve.
We are in the composite $N$ case so there is already a cyclic group 
inside the reduced group. 
For every maximal cyclic subgroup of the reduced automorphism group
we change coordinates so that 
the model of our curve is of the form:
\[
y^2=f(x)=\sum_{\nu=0}^s a_{2\nu} x^{2\nu}.
\] 
If there is an elementary abelian group inside the reduced group modulo 
some prime $p$ then  the reduced group contains  $E \rtimes(\Z /m\Z)$, 
for some $m$ prime to $p$.
The Galois cover $\mathbb{P}^1_k \rightarrow
\mathbb{P}^1_k$  with group $G=\mathbb{Z}_p^n \rtimes \langle \sigma \rangle$
is ramified over $0,\infty$. The elementary 
abelian group $E$ fixes either $0$ or $\infty$.
If $E$ fixes $\infty$ then the arbitrary 
$\tau \in E$ is of the form $\tau:x \mapsto x+c(\tau)$, 
for some $c(\tau)\in k$.  In this case 
$c$ is a root of the polynomial $G_1(x):=f(x+c_i)-f(x)$. 
If $E$ fixes $0$ then we consider the 
model $y^2=f^*(x)$ of the hyperelliptic curve and now $E$ fixes 
$\infty$. In this case $c$ is a root of the polynomial  $G_2(x):=f^*(x+c_i)-f^*(x).$

The possible primes $p$, so that there is  an elementary abelian $E$ group 
in the reduced group modulo $p$, are the divisors of $2g+2$, if $E$ 
does not fix any of the roots of the right hand side of the 
defining equation of the hyperelliptic curve. 
If $E$ fixes such  a root then $p \mid 2g+1$, but then 
there are a lot of cyclic groups in the reduced group.
This is the case in $A(48,7)$ and in $A(50,5)$.

Looking at the degrees of modular hyperelliptic curves $X_0(N)$,
where $N$ is composite we obtain that only for $p=3$ the reduction might have 
some extra automorphisms. The modular hyperelliptic curves with degree divisible by $3$ are those with $N\in\{ 22, 26, 28, 50, 46\}$. On the generic fibre 
if $N\neq 28$ then the  reduced group on the 
generic fibre is isomorphic to $\Z/2\Z$ while for $N=28$ the 
reduced group on the generic fibre is isomorphic to $D_4$.
 
We compute the polynomials $G_1(N)$ and $G_2(N)$ for every $N \in
{22,26,28,50,46}$ and we arrive at Table \ref{table:G1(N),G2(N)}.

\begin{table} 
 \caption{Tables of $G_1(N),G_2(N)$}
\extrarowheight3pt
{
\label{table:G1(N),G2(N)} 
\begin{tabular}{c|c|l}
\hline\hline
$22$ & $G_1$ &
$
\left( 2\,{c}^{3}+2\,c \right) {x}^{3}+ \left( 2\,{c}^{3}+c \right) x
+{c}^{6}+2\,{c}^{4}+2\,{c}^{2}
$ \\ 
& $G_2$ &
$
0
$ \\
\hline
$26$ & $G1$ & 
$
\left( 2\,{c}^{3}+c \right) {x}^{3}+x{c}^{3}+{c}^{6}+{c}^{4}
$ \\
& $G_2$ &
$
{x}^{3}{c}^{3}+2\,xc+2\,{c}^{6}+{c}^{2}
$ \\
\hline
$28$ & $G_1$ &
$
\left( {c}^{3}+2\,c \right) {x}^{3}+ \left( 2\,{c}^{3}+c \right) x+2
\,{c}^{6}+2\,{c}^{4}+2\,{c}^{2}
$ \\
& $G_2$ & 
$
\left( {c}^{3}+2\,c \right) {x}^{3}+ \left( 2\,{c}^{3}+c \right) x+2
\,{c}^{6}+2\,{c}^{4}+2\,{c}^{2}
$ \\
\hline
$28b$ & $G_1$ &
$
\left( {c}^{3}+2\,c \right) {x}^{3}+ \left( 2\,{c}^{3}+c \right) x+2
\,{c}^{6}+2\,{c}^{4}+2\,{c}^{2}
$ \\
& $G_2$ & 
$
\left( {c}^{3}+2\,c \right) {x}^{3}+ \left( 2\,{c}^{3}+c \right) x+2
\,{c}^{6}+2\,{c}^{4}+2\,{c}^{2}
$ \\
\hline
$28c$ & $G_1$ &
$
2\,{x}^{3}{c}^{3}+{c}^{6}+{c}^{3}
$ \\
& $G_2$ & 
$
2\,{x}^{3}{c}^{3}+{c}^{6}+{c}^{3}
$ \\
\hline
$50$ & $G_1$ & 
$
\left( {c}^{3}+c \right) {x}^{3}+ \left( {c}^{3}+c \right) x+2\,{c}^{
6}+{c}^{4}+2\,{c}^{2}
$ \\
& $G_2$ &
$
\left( {c}^{3}+2\,c \right) {x}^{3}+ \left( 2\,{c}^{3}+2\,c \right) x
+2\,{c}^{6}+2\,{c}^{4}+{c}^{2}
$ \\
\hline
$46$& 
 $G_1$ & $0$ \\
\hline
$37$ & $G_2$ & $0$ \\
\hline\hline
\end{tabular}
}
\end{table}

From that table we see that the curves $X_0(22)_3,X_0(28)_3,X_0(50)_3,X_0(46)_3$ admit an extra  $\Z/3\Z$ automorphism group, while $X_0(26)_3$ does not. 

\subsection{The $X_0(46)$ curve.}
The curve $X_0(46)$ is a curve of genus $5$ and 
has reduced automorphism group in the generic fibre 
isomorphic to $\Z/2\Z =\langle \sigma_1 \rangle$.
We compute the following model of our curve, 
so that $\sigma_1$ acts like $\sigma_1:x \mapsto -x$.
\begin{eqnarray*}
y^2 &=&
{x}^{12}+{\frac { \left( -4896\,\sqrt {2}-6786 \right) {x}^{10}}{-2392
\,\sqrt {2}-3381}}+{\frac { \left( -3512\,\sqrt {2}-4891 \right) {x}^{
8}}{-2392\,\sqrt {2}-3381}}-\,{\frac {2652{x}^{6}}{-2392\,\sqrt {2}-
3381}}  \\ 
&+ &
{\frac { \left( 3512\,\sqrt {2}-4891 \right) {x}^{4}}{-2392\,
\sqrt {2}-3381}}+{\frac { \left( 4896\,\sqrt {2}-6786 \right) {x}^{2}}
{-2392\,\sqrt {2}-3381}}+{\frac {-3381+2392\,\sqrt {2}}{-2392\,\sqrt {
2}-3381}}.
\end{eqnarray*}
We compute the decomposition of 
the coefficients of the above polynomial
and we see that there is no prime $P$ of $\Z[\sqrt{2}]$ with 
$N(P) \neq 2,3,23$ so that the reduction of the curve at this 
prime has a cyclic group containing $\sigma_1$.

We compute the dihedral invariants of this curve. 
We know that  if the reduced group $\bar{G}$ contains  
$(\Z/2\Z)^2$ as a subgroup, then  $2^{g-1}u_1^2 -u_g^{g+1}=0$ \cite[th. 3.8]{GutSha}.
The primes $p$ such that  $2^{g-1}u_1^2 -u_g^{g+1}$  becomes zero modulo $p$
are possible primes where the reduced group can be large enough to 
contain  $(\Z/2\Z)^2$. We compute that 
\[
2^{\frac{g-1}{2}} u_1 -u_g^{\frac{g+1}{2}}=
\frac{2^{23} \cdot 3^{12} \cdot 337^2 \cdot 683^2}{23^6 \cdot 
(-147+104 \sqrt{2})^5(147 +104\sqrt{2})
},
\]
\[
2^{\frac{g-1}{2}} u_1 +u_g^{\frac{g+1}{2}}=
\frac{2^{10} \cdot 3^{12} \cdot 7^2 \cdot 
13^2 \cdot 281^2 \cdot 709^2
}{23^6 \cdot 
(-147+104 \sqrt{2})^5(147 +104\sqrt{2})
}
\]
Only the prime $p=3$ gives rise to an extra automorphism
modulo $p$. For this prime we have reduced group isomorphic 
to $\Z/2\Z \times \Z/2\Z$ and $A(46,3)\cong(\Z/2\Z)^3$.

We will now see if there are primes so that 
$\sigma_1$ is the involution of a dihedral group and there 
is an element $\tau$ of order $n$ so that $\sigma_1,\tau$ 
generate a dihedral group. 
We consider the polynomials 
$g(x)=\frac{(x+1)^{2s}}{f(-1)}f\left(\frac{1-x}{1+x}\right)=
f_{A}^*(x)$, where $A=\left(\begin{smallmatrix} -1 & 1 \\ 1 & 1
\end{smallmatrix}\right)$, 
so that $\sigma_1$ acts on the model $y^2=g(x)$ as an 
involution sending $x\mapsto 1/x$.
The polynomial $g(x)$ is given by 
\[g(x)=g_1(x) g_2(x),\]
where
\[g_1(x)=
\left(
{x}^{6}+5/2\,\sqrt {2}{x}^{5}+5\,{x}^{4}+{\frac {21}{4}}\,\sqrt {2}{x}
^{3}+5\,{x}^{2}+5/2\,\sqrt {2}x+1
\right)
\]
\[
g_2(x)=
\left(
{x}^{6}+5/2\,\sqrt {2}{x}^{5}+7\,{x}^{4}+{\frac {25}{4}}\,\sqrt {2}{x}
^{3}+7\,{x}^{2}+5/2\,\sqrt {2}x+1
\right).
\]
A model of a dihedral group acting on the rational 
function field is given by $\tau_{norm}:x \mapsto \frac{1}{x}$ and 
$\sigma_{norm}:x \mapsto \zeta x$, where $\zeta$ is a primitive 
$n$-th root of one. 
This model is conjugate by a M\"obius transformation 
$Q:x\mapsto \frac{ax+b}{cx+d}$ to 
any other dihedral action on the rational function field.
We have normalized so that $Q \tau_{norm} Q^{-1}=\tau_{norm}$.
This means that for $a,b,c,d$ we have $b=\lambda c$, $a=\lambda d$, 
$d=\lambda a$, $c=\lambda b$ for some non zero element $\lambda$.
This gives that $\lambda=\pm 1$ therefore $a=\pm b$ and $b=\pm c$. 
If there is a dihedral group acting on our curve modulo $p$
then there should be a transformation 
of $Q$ of the form $x\mapsto \frac{ax+b}{\pm b x +\pm a}$ 
such that the polynomial $g^*_Q(x)$ has zero 
the coefficients $a_1,a_7,a_{11}$ of the monomials $x,x^7,x^{11}$.
The coefficients $a_1,a_7,a_{11}$ are polynomials 
of $a,b$ and we can eliminate $a$ by using the resultant determinant. 
We compute:
\[
\mathrm{Resultant}_a(a_1,a_7)=
\frac{3^{12}\cdot7^2 \cdot 23^{10}\cdot 60876058793276893^2}{2^{28}} b^{144}.
\]
The reduced  curve can have more automorphisms only at the primes dividing 
 the numerator of 
the resultant.
By reducing the hyperelliptic curve modulo all that primes 
and studying its automorphism group we see that the 
only prime where the reduced automorphism group grows is $p=3$. 

\section{The prime $N$-case}
\label{sec4}
If $N$ is a prime number $N\neq 37$ so that $X_0(N)$ is hyperelliptic then on the generic fibre 
there is only one involution, the hyperelliptic involution. 
In order to determine the possible primes $p\nmid N$ so that $X_0(N)_p$ has automorphism group 
greater than $\Z/2\Z$ we proceed as follows:
Suppose that the curve admits the hyperelliptic model $y^2=f_N(x)$ where $f_N(x) \in \Z[x]$. 
We consider an arbitrary M\"obius transformation $\sigma$  given by $x \mapsto \frac{ax +b}{cx+d}$.
Then we consider the coefficients of the polynomial 
\begin{equation} \label{equation: 71}
 f_N(x)-f_N\left( \frac{ax +b}{cx+d} \right) (cx+d)^{\deg{f_N}}=\sum_{\nu=0}^{\deg f_N } a_i x^i.
\end{equation}
If $\sigma$ is an automorphism then all $a_i$ should be zero. We would like to 
find if the system $a_i=0$ has solutions modulo $p$.
We consider the ideal $I_r:=\left\langle a_i, i=1,\ldots,r \right\rangle \lhd \Z [a,b,c,d]$ where $r < \deg f_N$.
We compute a Gr\"obner basis for the ideal $I_r$ with respect of the lexicographical order 
$a<b<d<c$, and then we form the set $S$ of all basis elements that are polynomials in $c$ only. 
Since in the generic fibre the only admissible automorphism is the trivial one, 
the greatest common divisor of elements in $S$ is $c^\alpha$ for some $1<\alpha \in \mathbb{N}$.
We divide every element in $S$ by $c^\alpha$ and we obtain an integer $\delta$ as an element in 
the set $\{f/c^\alpha : f \in S\}$. The prime factors $p$ of $\delta$ are exactly the possible primes 
where an automorphism $\sigma $ with $c\neq 0$ can appear. 

We do again the same procedure, but now we choose the a lexicographical order where $a<c<d<b$. 
We find again an integer $\delta'$ and the divisors of $\delta'$ are exactly the 
possible primes where an automorphism $\sigma$ with $d\neq 0$ can appear.

If we select a big $r$ then the procedure of finding the Gr\"obner 
basis is difficult. If on the other hand we select a small $r$ ($4 \leq r$ since we need 
at least four equations in order to find a unique solution in $a,b,c,d$) then 
we  obtain  big integers $\delta,\delta'$ that we are 
not able to factorize. A selection of $r=6$ allows us to perform the computations needed for the 
given set of $N$. 

Now for each prime $p$ that is a divisor of $\delta$ or $\delta'$ we consider 
the ideals $I_{\deg f_N} \otimes_\Z \Z/p\Z$ and we do the same elimination 
procedure. The Gr\"obner basis computation is easier to perform over a finite field 
and we finally arrive at a solution of the system $a_i =0 \mod p$.

For example, for the $N=41$ case the only exceptions can happen 
at the primes $2,17,41$. The primes $2,41$ are excluded so we focus to 
the $p=17$ case. We reduce our curve modulo $17$ and then we compute that the ideal $I_{\deg f_{41}} \otimes_\Z \Z/p\Z$ 
has a Gr\"obner basis of the form:
\[
 \{a+16d+b,d^8+12b^8+16,b(d+8b),c+8b,b(b^8+13) \}.
\]
We will now solve the above system.
If $b=0$ then we see that $c=0$ and $a=d$,  therefore we obtain the 
identity matrix.
If $b\neq 0$ then $b^8+13=0 \Rightarrow b^4=2$.
Let $b$ be a fourth root of $2$ in $\bar{\mathbb{F}}_{17}$.
Then $c=-8b$, $d=-8b$, $a=-9b$. The equation 
$d^8+12 b^8+16$ is compatible with the system.
Thus we obtain the extra automorphism $\sigma$ so that 
$\bar{\sigma}:x \mapsto (-9b x+b)/(-8bx-9b)=(9x-1)/(8x+9)$.
The automorphism group in this case is $(\Z/2\Z) \times (\Z/2\Z)$.

There might be also an extra automorphism $\sigma$ modulo $p$ so that 
$b=c=0$. Then $\sigma: x \mapsto(a/d)x$ where $a/d$ is an $n$-th root of unity. 
This does not happen for any prime $p$ as one easily checks.  
Our results are collected  in table \ref{Table: Results}.

\begin{rem} For the case $N=71$, the polynomial $f_{71}(x)$ has degree
  $14$. Then the coefficients $a_i$ in \eqref{equation: 71} are
  polynomials of degree $14$ in $a,b,c,d$. The computation of the
  Gr\"obner basis over $\Z$ is very time and memory-consuming. For
  such situations, we use the following trick. If we homogenize
  the polynomial $f_N(x)$ into a binary form $f_N(x,y)$, then the
  property \eqref{equation: 71} means that $f_N(x,y)$ is invariant
  under the substitution $(x,y)\mapsto(ax+by,cx+dy)$. Any
  transvectant
  $$
    (f_N,f_N)^r=\sum_{k=0}^r(-1)^k\binom rk
    \frac{\partial^r f_N}{\partial x^{r-k}\partial y^k}
    \frac{\partial^r f_N}{\partial x^k\partial y^{r-k}}
  $$
  of $f_N(x,y)$ with itself is also invariant under this
  substitution. (See \cite[page 54]{Glenn}.) Thus, for $X_0(71)$ we
  first work on $g=(f_{71},f_{71})^{10}$ of degree $8$. We determine
  all possible primes $p$ such that there exists non-trivial
  automorphisms for $g$ modulo $p$. We then determine whether these
  primes indeed give extra automorphisms for $X_0(71)_p$.
\end{rem}

%

%
%

\section{Automorphisms in characteristic $2$}
\label{section: char 2}

In this section we will study modular hyperelliptic curves $X_0(N)$ for $N$ 
odd, in characteristic $2$. These curves admit a minimal Weierstrass
model of the form $y^2+q(x)y=p(x)$ \cite{Lockhart}, that can easily be
found with the help of {\em magma} algebra system.  The equations are
given in Table \ref{table: mod 2 equations}. These models are examples 
of the unified Artin-Schreier-Kummer theory in the sense of 
T. Sekiguchi, N. Suwa \cite{SeSu94}, \cite{SeSu95}.

\begin{table} \label{table: mod 2 equations}
\extrarowheight3pt
\caption{Global minimal Weierstrass equations for $X_0(N)$, $N$ odd}
{\tiny
\begin{tabular}{|c||l|} \hline\hline
$N$ & Equation \\ \hline\hline
$23$ & $y^2 + (x^3+ x + 1)y = -2x^5 - 3x^2 + 2x - 2$ \\ \hline
$29$ & $y^2 + (x^3 + 1)y = -x^5 - 3x^4 + 2x^2 + 2x - 2$ \\ \hline
$31$ & $y^2 + (x^3 + x + 1)y = -2x^5 + x^4 + 4x^3 - 
    3x^2 - 4x - 1$ \\ \hline
$33$ & $y^2 + (x^4 + x^2 + 1)y = 2x^7 + 9x^6 + 27x^5
    + 56x^4 + 81x^3 + 85x^2 + 54x + 20$ \\ \hline
$35$ & $y^2 + (x^4 + x^2 + 1)y = -x^7 - 2x^6 - x^5 - 
    3x^4 + x^3 - 2x^2 + x$ \\ \hline
$37$ & $y^2 + (x^3 + x^2 + x + 1)y = 3x^5 + 8x^4 + 
    11x^3 + 8x^2 + 3x$ \\ \hline
$39$ & $y^2 + (x^4 + x^3+ x^2 + x + 1)y = -2x^7 + 
    2x^5 - 7x^4 + 2x^3 - 2x$ \\ \hline
$41$ & $y^2 + (x^4 + x)y = -x^7 - 2x^6 + 2x^5 + 5x^4
    + 2x^3 - 4x^2 - 5x - 2$ \\ \hline
$47$ & $y^2 + (x^5 + x^3 + 1)y = x^9 - 8x^7 - 34x^6 -
    74x^5 - 106x^4 - 103x^3 - 67x^2 - 25x - 5$ \\ \hline
$59$ & $y^2 + (x^6 + 1)y = x^{11} - 7x^9 - 21x^8 - 
    38x^7 - 51x^6 - 53x^5 - 44x^4 - 30x^3 - 17x^2 - 6x - 3$ \\ \hline
$71$ & $y^2 + (x^7 + x^6 + x^4 + x + 1)y = -3x^{13} + 
    9x^{12} - 17x^{11} + 16x^{10} - 12x^9 +$ \\
 & $+3x^8 + 9x^7 - 17x^6 + 16x^5 - 
    15x^4 + 10x^3 - 7x^2 + x - 2$\\ \hline\hline
\end{tabular}
}
\end{table}

\begin{lemma} \label{lemma: automorphism over 2} Let $\mathscr
  C:=y^2+q(x)y+p(x)$ be a hyperelliptic curve of genus $g$ over
  $\overline\F_2$ with $\deg q(x)\le g+1$ and $\deg p(x)\le
  2g+1$. Then every automorphism $\sigma$ of $\mathscr C$ is of the form
  $$
    \sigma:(x,y)\longmapsto\left(
    \frac{ax+b}{cx+d},\frac{y+h(x)}{(cx+d)^{g+1}}\right)
  $$
  for some $\left(\begin{smallmatrix}a&b\\c&d\end{smallmatrix}\right)
  \in\GL(2,\overline\F_2)$ and $h(x)\in\overline\F_2[x]$ of degree at
  most $g+1$ satisfying
  $$
    q\left(\frac{ax+b}{cx+d}\right)(cx+d)^{g+1}=q(x), \quad
    p\left(\frac{ax+b}{cx+d}\right)(cx+d)^{2g+2}=p(x)+h(x)^2+q(x)h(x).
  $$
  In particular the hyperelliptic involution is given by 
\begin{equation} \label{hyp-inv}
j(x)=x, \quad j(y)=y+q(x).
\end{equation}
\end{lemma}

\begin{proof} The function field of $\mathscr C$ is an Artin-Schreier
  extension of the rational function field $\overline\F_2(x)$. Indeed,
  if we set $Y=y/q$, then we have
  $$
    Y^2+Y=\frac p{q^2},
  $$
  and the hyperelliptic involution is given by $(x,Y)\mapsto(x,Y+1)$, i.e,
  $\sigma(y)=y+q$.
  The hyperelliptic involution is in the center of the automorphism
  group. Thus, the restriction of an automorphism $\sigma$ of
  $\mathscr C$ to $\overline\F_2(x)$ gives an automorphism of
  $\overline\F_2(x)$. Therefore, we must have $\sigma(x)=(ax+b)/(cx+d)$
  for some $\left(\begin{smallmatrix}a&b\\c&d\end{smallmatrix}\right)$
  in $\GL(2,\overline\F_2)$.

  Recall \cite[prop. 1.12]{Lockhart} that a basis for the space of holomorphic differentials on
  $\mathscr C$ is given by
  $$
    \omega_i=\frac{x^{i-1} dx }{2y + q}=\frac{x^{i-1}dx}{q}, \;\;
    1 \leq i  \leq g,
  $$
  and every automorphism $\sigma$ of $\mathscr C$ induces a linear
  action on the space of holomorphic differentials. Write
  $q((ax+b)/(cx+d))(cx+d)^{g+1}=q^\ast(x)\in\overline\F_2[x]$. We find
  $$
    \sigma(\omega_i)=\frac{\sigma(x)^{i-1}d\sigma(x)}{q(\sigma(x))}
   =(ad-bc)(ax+b)^{i-1}(cx+d)^{g-i}\frac{dx}{q^\ast}.
  $$
  Since each $\sigma(\omega_i)$ is a linear combination of $\omega_j$,
  we must have $q^\ast=\lambda q$ for some
  $\lambda\in\overline\F_2^\ast$. Because for any
  $\alpha\in\overline\F_2^\ast$,
  $\alpha\left(\begin{smallmatrix}a&b\\c&d\end{smallmatrix}\right)$
  defines the same automorphism on $\overline\F_2(x)$ as
  $\left(\begin{smallmatrix}a&b\\c&d\end{smallmatrix}\right)$, we may
  rescale $a,b,c,d$ so that $\lambda=1$, i.e., we have
  \begin{equation} \label{equation: q^ast}
    q\left(\frac{ax+b}{cx+d}\right)(cx+d)^{g+1}=q(x).
  \end{equation}

  We now consider $\sigma(y)$. We write it in the form $\sigma(y)=\mu
  y+\nu$ with $\mu,\nu\in\overline\F_2(x)$. Substituting the
  expression into $\sigma(y)^2+q(\sigma(x))\sigma(y)+p(\sigma(x))=0$ and
  using \eqref{equation: q^ast}, we obtain
  $$
    y^2+\frac{q(x)}{\mu(cx+d)^{g+1}}y+\frac{p(\sigma(x))+\nu^2
   +q(x)\nu/(cx+d)^{g+1}}{\mu^2}=0.
  $$
  Comparing this with $y^2+q(x)y+p(x)=0$, we find $\mu=1/(cx+d)^{g+1}$
  and $\nu(cx+d)^{g+1}$ is a polynomial $h(x)$ such that
  $$
    p\left(\frac{ax+b}{cx+d}\right)(cx+d)^{2g+2}
   =p(x)+h(x)^2+q(x)h(x).
  $$
  This completes the proof of the lemma.
\end{proof}

\begin{lemma}
Let $X$ be a hyperelliptic curve in characteristic $2$.
The group structure of the full automorphism subgroup $G$  of $X$
is determined by the structure of the $2$-Sylow subgroup of $G$.
\end{lemma}
\begin{proof}
Let $\bar{G}$ denote the reduced group of $G$.
By the theory of group extensions, the group $G$ is determined uniquely by a cohomology 
class in the group $H^2(\bar{G},\Z/2\Z)$ corresponding to the first row of 
eq. \ref{subext}.

For $p$-prime let $H^2(\bar{G},\Z/2\Z)_p$ denote the $p$-part of the 
finite abelian group $H^2(\bar{G},\Z/2\Z)$ and let $\bar{G}_p$ 
denote the $p$-Sylow subgroup of the reduced group
$\bar{G}$.
The following map is monomorphism: 
\[
\Phi: H^2(\bar{G},\Z/2\Z) =\bigoplus_{p |s} H^2(\bar{G},\Z/2\Z)_p
\rightarrow \bigoplus_{p |s } H^2(\bar{G}_p,\Z/2\Z),
\]
\[
\alpha=\sum_{p|s} \alpha_p \mapsto 
\sum_{p|s} res_{\bar{G},\bar{G}_p} \alpha_p
\]
\cite[p. 93]{WeissCoho}.
Now if $(p,2)=1$ then $H^2(\bar{G}_p,\Z/2\Z)=0$.
Therefore we have to consider only the 
$p=2$ case. This proves that $H^2(\bar{G},\Z/2\Z)$ is 
itself a group of order 2.
Moreover, the  restriction map 
$H^2(\bar{G},\Z/2\Z) \rightarrow  H^2(\bar{G}_2,\Z/2\Z)$
is a monomorphism.
The class $res_{\bar{G},\bar{G}_2}(\alpha)$ corresponds to the subextension 
given by the second row of equation \ref{subext}
\begin{equation} \label{subext}
\xymatrix{
1 \ar[r] &  \Z/2\Z \ar[r] \ar_{=}[d] & G  \ar^{\pi}[r]   & \bar{G} \ar[r] & 1 \\
1 \ar[r] &  \Z/2\Z \ar[r] & G_2 \ar[r] \ar[u]_{1-1} & \bar{G}_2 \ar[r] \ar[u]_{1-1} & 1.
}
\end{equation}
This means that the structure of $G_2$ determines uniquely the structure of $G$.

%
%
%

\end{proof}

\begin{cor}
If $\bar{G}_2=\{1\}$ then $G=\Z/2\Z \times \bar{G}$.
\end{cor}

\begin{cor} \label{corram}
Let $\bar{G}_2$ denote the $2$-Sylow subgroup of the reduced group. This group 
is elementary abelian and  fixes only one point  $\infty$ of $\mathbb{P}^1$.
If $\infty$ does not ramify in $X \rightarrow \mathbb{P}^1$ then 
$G=\bar{G} \times \Z/2\Z$.
\end{cor}
\begin{proof}
If $\infty$ does not ramify in $X \rightarrow \mathbb{P}^1$ then 
there are two points  $P_1,P_2$ of $X$  above $\infty$.
The group $\langle \sigma \rangle:=\Z/2\Z=\mathrm{Gal}(X/\mathbb{P}^1)$ transfers 
$P_1$ to $P_2$.
Consider the isotropy supgroup $G_2(P_1)$. By the 
transitivity of the ramification index we have that $|G_2(P_1)|=|\bar{G}_2|$.
The map $G_2(P_1) \rightarrow \bar{G}_2$ is onto and since the two groups
have the same order it is an isomorphism. Therefore, $\pi$ has a section 
and last line of (\ref{subext}) splits. The first line splits also and 
$G=\bar{G} \times \Z/2\Z$, since the hyperelliptic involution is central. 
\end{proof}

The group $\bar{G}_2$ is elementary abelian therefore is it is generated by the 
commuting elements $x_i$, $i=1,\ldots,s$.
If the group $G_2$ is not a direct product then there are elements say $\sigma \in G_2$ 
of order $4$. 
Every such element when raised to the square gives the hyperelliptic involution, 
i.e. 
 $j=\sigma^2$. 
For the elements $x_i$ such that $\pi^{-1}\langle x_i \rangle$ is a cyclic 
group of order $4$ select a generator $\sigma_i$ for this group.
All these elements have squares equal to $j$.
Order the elements $x_i$ so that for $1\leq i \leq \nu_0$  the group  $\pi^{-1}\langle x_i \rangle$ is a cyclic 
group of order $4$ and for $\nu_0<1$ the group $\pi^{-1}(\langle x_i \rangle)$ is a direct product of two 
cyclic groups. 
The group $G_2$ admits the following presentation 
in terms of generators and relations:
\[
G_2:=\left\langle j, \sigma_\nu , x_\mu : 1 \leq \nu \leq \nu_0 < \mu \left| 
\begin{array}{l}
\sigma_\nu^2=j,j^2=1,x_\mu^2=1, \\ 1=[x_\mu,x_{\mu'}]=[x_\mu,\sigma_\nu^2]=[\sigma_\nu,\sigma_{\nu'}]
 \end{array}
 \right.
 \right\rangle.
\]

For every Weierstrass model given in \ref{table: mod 2 equations} we use 
lemma \ref{lemma: automorphism over 2} in order to determine 
the automorphism group in characteristic $2$. For every entry in that 
table we do not get any new automorphism except in the cases $N=33,37$.

\subsection{Case $X_0(37)$} By Table \ref{table: mod 2 equations},
a Weierstrass model for $X_0(37)_{/2}$ is given by
$y^2+q(x)y=p(x)$ with
$$
  q(x)=x^3+x^2+x+1,\qquad
  p(x)=x^5 + x^3 + x.
$$
The hyperelliptic involution $j$ is $j:(x,y)\mapsto(x,y+q(x))$. Let
$G$ denote the automorphism group of $X_0(37)_{/2}$ and $\overline
G=G/\gen j$ be the reduced automorphism group, considered as a
subgroup of $\Aut\overline\F_2(x)=\operatorname{PGL}(2,\overline\F_2)$.
According to Lemma \ref{lemma: automorphism over 2}, an
element of $G$ takes the form
$$
  (x,y)\longmapsto\left(\frac{ax+b}{cx+d},\frac{y+h(x)}{(cx+d)^3}\right),
$$
where $h(x)=u_0+u_1x+u_2x^2+u_3x^3\in\overline\F_2[x]$ is a polynomial
of degree at most $3$ and
$\left(\begin{smallmatrix}a&b\\c&d\end{smallmatrix}\right)\in
\overline G\subset\GL(2,\overline\F_2)$ satisfies
$$
  q\left(\frac{ax+b}{cx+d}\right)(cx+d)^3=q(x), \quad
  p\left(\frac{ax+b}{cx+d}\right)(cx+d)^6=p(x)+h(x)^2+q(x)h(x).
$$
These two conditions give a set of relations among $a,b,c,d$ and
$u_i$. The Gr\"obner basis of the ideal generated by these relations
with respect to the lexicographic order $u_0>u_1>u_2>u_3>a>b>d>c$
is
\begin{equation*}
\begin{split}
  &u_0 + u_3 + d^2c^4 + d^2c + dc^8 + dc^2 + c^{192} + c^{180} + c^{168}
    + c^{165} + c^{150} + c^{138} + c^{135} \\
  &\quad + c^{132} + c^{120} + c^{105} + c^{96} + c^{90} + c^{84} + 
        c^{75} + c^{69} + c^{66} + c^{48} + c^{36} + c^{18} + c^9, \\
  &u_1 + u_3 + d^2c + dc^8 + c^{168} + c^{138} + c^{120} + c^{105} +
    c^{90} + c^{75} + c^{72} + c^{60} + c^{48} \\
  &\quad + c^{45} + c^{30} + c^{24} + c^{18} + c^{15} + c^{12} + c^3,
  \\
  &u_2 + u_3 + d^2c^4 + dc^2 + c^{180} + c^{165} + c^{150} + c^{144}
    + c^{135} + c^{129} + c^{96} + c^{84} + c^{69} \\
  &\quad + c^{66} + c^{60} + c^{48} + c^{45} + c^{36} + c^{33} + c^{30} + 
        c^{18} + c^{15}, \\
  &u_3^2 + u_3 + d^2c^4 + d^2c + dc^5 + dc^2 + c^{36} + c^{33} +
    c^{21} + c^{18} + c^6 + c^3, \\
  &a + d + c^{16} + c, \\
  &b + c^{16}, \\
  &d^3 + d^2c + dc^2 + c^{192} + c^{144} + c^{132} + c^{129} + c^{72}
     + c^{48} + c^{33} + c^{24} + c^{18} \\
  &\quad + c^{12} + c^9 + 1, \\
  &d(c^{16} + c) + c^{176} + c^{161} + c^{146} + c^{131} + c^{80} +
    c^{65} + c^{56} + c^{41} + c^{26} + c^{20} \\
  &\quad + c^{17} + c^{11} + c^5 + c^2, \\
  &(c^{16}+c)(c^{192}+c^{144}+c^{132}+c^{129}+c^{96}+c^{72}+c^{66}+c^{48}
    +c^{36}+c^{33} \\
  &\quad+c^{24}+c^{18}+c^{12}+c^9+c^6+c^3+1).
\end{split}
\end{equation*}
Here the polynomial of degree $192$ in $c$ in the last element of the
Gr\"obner basis is a product of $12$ irreducible polynomials of degree
$8$ over $\F_2$. Using this basis, we find that the total number of
solutions in $\overline \F_2$ is $480$. (Each root of the degree $192$
polynomial gives two solutions and each root of $c^{16}+c$ gives $6$
solutions.) However, since for each root $\alpha$ of $x^3+1$ in
$\F_4$, $(u_0,u_1,u_2,u_3,a,b,c,d)$ and $(u_0,u_1,u_2,u_3,\alpha
a,\alpha b,\alpha c,\alpha d)$ give the same automorphism, we find
that
$$
  |G|=480/3=160, \qquad |\overline G|=|G|/2=80.
$$
We now determine the structure of the automorphism group.

We first consider the reduced automorphism group $\overline G$.
Recall that, in general, the order of a matrix in
$\operatorname{PGL}(2,\overline\F_2)$ can only be $2$ or an odd
integer. Moreover, the order is $1$ or $2$ if and only if the trace
is zero. Now the relation $a+d+c^{16}+c=0$ shows that an
element $\left(\begin{smallmatrix}a&b\\c&d\end{smallmatrix}\right)
\in\overline G$ has order $1$ or $2$ if and only if $c\in\F_{16}$.
Therefore, we find that the Sylow $2$-subgroup of $\overline G$ is
an elementary abelian $2$-group of order $16$, and is normal in
$\overline G$. (Again, each $c\in\F_{16}$ gives $3$ solutions
$(a,b,c,d)$, but for each root $\alpha$ of $x^3+1$ in $\overline\F_2$,
$(a,b,c,d)$ and $(\alpha a,\alpha b,\alpha c,\alpha d)$ correspond to
the same reduced automorphism in $\overline G$.) The remaining
$64$ elements of $\overline G$ all have order $5$, and $\overline G$
has $16$ Sylow $5$-subgroups. Therefore,
$\overline G$ is the semi-direct product of an elementary abelian
$2$-group of order $16$ by a cyclic group of order $5$.

Now consider the structure of $G$ itself. Let $P$ be its Sylow
$2$-subgroup, and $\tau$ be any element of order $5$. The centralizer
$Z_P(\tau)$ of $\tau$ in $P$ must satisfy
$|Z_P(\tau)|\equiv|P|=32\mod 5$. Thus, we have $|Z_P(\tau)|=2$ or
$|Z_P(\tau)|=32$. The latter possibility cannot occur as it would imply
that $\overline G$ is an abelian group. Thus, we have $|Z_P(\tau)|=2$,
that is, $Z_P(\tau)=\gen j$.

We next turn the attention to the center $Z(P)$ of the Sylow
$2$-subgroup $P$ itself. Observe that $\gen\tau$ acts on $Z(P)$ by
conjugation. The identity automorphism and the hyperelliptic
involution are left fixed by this group action. Since $|Z_P(\tau)|=2$,
all the other orbits under this group action have $5$ elements. In
other words, we have $|Z(P)|\equiv 2\mod 5$, i.e., $|Z(P)|=2$ or
$|Z(P)|=32$.

Assume that $|Z(P)|=32$. Then $P$ is an abelian group of order $32$
whose elements have order at most $4$ (since $P/\gen j$ is elementary
abelian). Noticing that $\tau$ acts on the set of elements of order $4$
in $P$ and that $|Z_P(\tau)|=2$, we see that the number of elements of
order $4$ in $P$ must be a multiple of $5$. The only possibility is
that $P\simeq(\Z/2\Z)^5$. However, we can easily check that the
automorphism $\sigma\in G$ given by
$$
  \sigma:(x,y)\longmapsto\left(\frac{(\alpha+1)x+\alpha}
  {\alpha x+(\alpha+1)},\frac{y+x^2+x}{(\alpha x+(\alpha+1))^3}\right)
$$
is an element of order $4$, where $\alpha$ is a root of $x^2+x+1$ in
$\F_4$. Therefore, we conclude that $|Z(P)|$ cannot be $32$. Instead,
we have $|Z(P)|=2$, i.e., $Z(P)=\gen j$.

Now we have $|Z(P)|=2$ and $P/Z(P)$ is elementary abelian. This means
that $P$ is one of the extraspecial groups. (An \emph{extraspecial
  group} $H$ is a $p$-group such that the center $Z$ is cyclic of
order $p$ and the quotient group $H/Z$ is a non-trivial elementary
abelian $p$-group.) For order $32$, there are two extraspecial groups
$E_{32+}:=(D_4\times D_4)/\gen{(a,a)}$ and
$E_{32-}:=(D_4\times Q_8)/\gen{(a,b)}$, where $a$ and $b$ denote the
non-trivial elements in the centers of the dihedral group $D_4$ and
the quaternion group $Q_8$, respectively. To determine which one $P$
is isomorphic to, we consider the action of $\gen\tau$ defined by
conjugation on the set $S$ of subgroups of order $16$ in $P$.

We claim that $\tau A\tau^{-1}\neq A$ for all $A\in S$. Assume that
$\tau A\tau^{-1}=A$ for some $A\in S$. The centralizer $Z_A(\tau)$ of
$\tau$ in $A$ must have only one element since $|Z_P(\tau)|=2$ and
$|Z_A(\tau)|\equiv|A|\mod 5$. In other words, $j\not\in A$, but this
would imply that $P\simeq A\times\gen j$, which cannot be true for an
extraspecial group. Therefore, we must have $\tau A\tau^{-1}\neq A$.
It follows that for each group $H$ of order $16$, the number of
subgroups of $P$ that are isomorphic to $H$ must be divisible by $5$.
Now according to the database of small groups \cite{Green}, $E_{32+}$
has $9$ subgroups isomorphic to $D_4\times(\Z/2\Z)$ and $6$ subgroups
isomorphic to another group $H_{16}$ of order $16$. Therefore, we conclude
that $P$ must be isomorphic to $E_{32-}$, which has $5$ subgroups
isomorphic to $Q_8\times(\Z/2\Z)$ and another $5$ subgroup isomorphic
to $H_{16}$, and the automorphism group $G$ is a semi-direct product
of $E_{32-}$ by a cyclic group of order $5$.

Of course, the conclusion above can be verified by brute force
computation. However, the computation is too complicated to be
presented here.

\subsection{Case $X_0(33)$}
The conditions given in lemma \ref{lemma: automorphism over 2} give rise to a system in $a,b,c,d$ describing 
every element in the reduced group $\bar{G}$. The Gr\"obner basis 
of this system is computed to be:
\[
 a^2 + ac + d^2 + dc,                                                                                                                                 
    ab + ac + bd + dc,                                                                                                                                 
    ad + bc + d^2 + dc + c^2,                                                                                                                           
    ac^4 + a + dc^4 + d,\]\[                                                                                                                                 
    b^2 + bd + dc + c^2,                                                                                                                                 
    bdc + bc^2 + dc^2 + c^3,                                                                                                                           
    d^4 + dc^3 + c^4 + 1,                                                                                                                                 
    d^2c + dc^2,                                                                                                                                         
    c^5 + c     
\]
This gives us the following solutions (written in matrix form):
\[
\bar{G}= \left\{
\begin{pmatrix}
 1 & 0 \\
0 & 1
\end{pmatrix},
\begin{pmatrix}
 0 &1 \\
1 & 0
\end{pmatrix},
\begin{pmatrix}
 0 & 1 \\
1 & 1
\end{pmatrix},
\begin{pmatrix}
 1 & 1 \\
0 & 1
\end{pmatrix},
\begin{pmatrix}
 1 & 1 \\
1 & 0
\end{pmatrix},
\begin{pmatrix}
 1 & 0 \\
1 & 1
\end{pmatrix}
\right\}.
\]
Thus, the group $\bar{G}$ is isomorphic to the group $\mathrm{GL}(2,\mathbb{F}_2)$ of order $6$.
One 2-Sylow subgroup of $\bar{G}$ is given by the group generated by  the element $\tau:x \mapsto x+1$.
The fixed point of $\tau$ is the point $\infty$  and since $\infty$ is not ramified in the 
cover $X \rightarrow X^{\langle j \rangle}$ Corollary \ref{corram} implies that $G=\bar{G} \times \langle j \rangle$.

\bibliographystyle{amsplain}
\def\cprime{$'$}
\providecommand{\bysame}{\leavevmode\hbox to3em{\hrulefill}\thinspace}
\providecommand{\MR}{\relax\ifhmode\unskip\space\fi MR }
\providecommand{\MRhref}[2]{%
  \href{http://www.ams.org/mathscinet-getitem?mr=#1}{#2}
}
\providecommand{\href}[2]{#2}

\end{document}